# A transient Markov chain with finitely many cutpoints

Nicholas James[1], Russell Lyons[*,2] and Yuval Peres[†,3]

*Berkeley, California, Indiana University and University of California, Berkeley*

**Dedicated to David Freedman with admiration**

**Abstract:** We give an example of a transient reversible Markov chain that almost surely has only a finite number of cutpoints. We explain how this is relevant to a conjecture of Diaconis and Freedman and a question of Kaimanovich. We also answer Kaimanovich's question when the Markov chain is a nearest-neighbor random walk on a tree.

## 1. Introduction

While studying extensions of De Finetti's theorem to Markov chains, Diaconis and Freedman [3] stated a general conjecture for transient Markov chains $\{S_n\}$. We give a result on cutpoints that is relevant to their conjecture. We begin with some background.

We say that an event $A$ in the space of trajectories of the Markov chain is **exchangeable** if it is invariant under finite permutations, i.e., if $(S_0, S_1, \ldots) \in A$, then so is $(S_{\pi(0)}, \ldots, S_{\pi(n)}, S_{n+1}, \ldots)$ for any $n$ and any permutation $\pi$ of $\{0, \ldots, n\}$. The $\sigma$-field of exchangeable events, $\mathcal{E}$, is called the exchangeable $\sigma$-field. Let $\overline{\mathcal{E}}$ be the completion of $\mathcal{E}$. A transient process visits each state only finitely often, and so for each state $x$ in the state space $X$ there is a random variable $V(x)$ that counts the number of visits, $V(x) := \#\{n \geq 0\,;\ S_n = x\}$. We call the collection $V := \{V(x)\}_{x \in X}$ the **occupation numbers** of the process. Clearly, $V$ is $\mathcal{E}$-measurable. A natural question, posed by Kaimanovich [6], is to determine under what conditions the exchangeable $\sigma$-field is generated by $V$. This was motivated by similar issues arising in the study [7] of random walks on lamplighter groups.

Write $V_n(x) := \#\{k \in [0, n]\,;\ S_k = x\}$. Note that an event $A \in \sigma(S_j\,;\ j \geq 0)$ is invariant under permutations of $S_0, \ldots, S_n$ if and only if $A \in \sigma(V_n, S_{n+1}, S_{n+2}, \ldots)$. Therefore

$$(1) \qquad \mathcal{E} = \bigcap_n \sigma(V_n, S_{n+1}, S_{n+2}, \ldots).$$

[*]Supported in part by NSF Grant DMS-04-06017.
[†]Supported in part by NSF Grant DMS-06-05166.
[1]Berkeley, CA, USA, e-mail: nicholasjames@cal.berkeley.edu
[2]Indiana University, Department of Mathematics, Bloomington, IN 47405-5701, USA, e-mail: rdlyons@indiana.edu
[3]University of California, Departments of Statistics and Mathematics, Berkeley, CA, 94720-3860 and Microsoft Corporation, One Microsoft Way, Redmond, WA 98052-6399, USA, e-mail: peres@stat.Berkeley.EDU

*AMS 2000 subject classifications:* Primary 60J10; secondary 60J50.

*Keywords and phrases:* birth-and-death chain, cutpoints, exchangeable, nearest-neighbor random walk, occupation numbers, transient Markov chain, trees.

For any Markov chain $\{S_n\}$, the sequence of transitions $\{(S_n, S_{n+1})\}$ is also a Markov chain; for such chains of transitions, Kaimanovich's question was posed earlier as a conjecture by Diaconis and Freedman in [3]. To be precise, let $M_n(x, y)$ be the number of transitions made from $x$ to $y$ up to time $n$, so that $M_n(x, y)$ increases to a finite limit $M(x, y)$ as $n \to \infty$. They made the following conjecture in [3]:

**Conjecture 1.1.** *The intersection of the $\sigma$-fields*

$$\bigcap_n \sigma(M_n, S_{n+1}, S_{n+2}, \ldots) \tag{2}$$

*is always generated (up to completion) by $M$.*

By comparing (2) to (1), we see that (2) is just the exchangeable $\sigma$-field for the chain of transitions $\{(S_n, S_{n+1})\}$.

James and Peres [5] related the questions above to *cutpoints* of the Markov chain trajectory. Call $x$ a **cutpoint** if for some $k$, we have $S_k = x$ and the future of the chain, $\{S_{k+1}, S_{k+2}, \ldots\}$, is disjoint from its past $\{S_0, S_1, \ldots, S_k\}$. Call $S_k$ a **strong cutpoint** if the probability of a transition from $S_i$ to $S_j$ is 0 whenever $i < k < j$. In [5], Conjecture 1.1 was proved under the condition

(3)   the Markov chain $\{S_n\}$ has infinitely many cutpoints almost surely.

We give a brief outline to illustrate the connection; see [5] for more details. Under the assumption (3), the portions $\psi_1, \psi_2, \ldots$ of the space-time path $(n, S_n)$ between successive cutpoints are conditionally independent given $M$, and the intersection (2) is contained in the tail $\sigma$-field of the $\{\psi_j\}_{j \geq 1}$, which is trivial (given $M$) by Kolmogorov's zero-one law. Conditional triviality of a $\sigma$-field given $M$ means that the $\sigma$-field is generated by $M$ up to completion.

James and Peres [5] also showed that if $\{S_n\}$ almost surely has infinitely many strong cutpoints, then $\overline{\mathcal{E}}$ is generated by the occupation numbers. Thus, if every transient Markov chain had infinitely many strong cutpoints a.s., then Kaimanovich's question would be resolved.

In general, one expects that a random walk that is "very transient" will have infinitely many strong cutpoints. As shown in [1, 5, 8], transient random walks on Cayley graphs have infinitely many strong cutpoints a.s. More precisely, Lawler [8] proved (3) for simple random walk on the lattices $\mathbf{Z}^d$ for $d \geq 4$ and his argument applies to strong cutpoints and to any Cayley graph with volume growth at least polynomial of degree 5. This was extended, using a different argument, to $\mathbf{Z}^3$ in [5]. Blachère [1] extended the argument of [5] and showed that simple random walks on all transient Cayley graphs of groups have infinitely many strong cutpoints.

This raises the natural question of whether *every* transient Markov chain has infinitely many cutpoints a.s.; a positive answer would establish the conjecture of Diaconis and Freedman. In Section 3 we show, however, that this is not true, even for birth-and-death chains.

## 2. Exchangeability, transition counts and trees

In this section, we show that for transient nearest-neighbor walks on trees, the exchangeable $\sigma$-field is generated by the occupation numbers. This result was established in the thesis [4] of the first author, but was never published; the proof here is shorter than in [4], but relies on the same ideas. Note that the example in Section 3 is a nearest-neighbor random walk on a special tree (a halfline) such that the walk a.s. has finitely many cutpoints, so the proof cannot rely on cutpoints.

Consider a transient Markov chain as in the introduction. If $V(x) > 0$, let $U(x)$ be the state visited by the Markov chain immediately after its last visit to $x$. For completeness, define $U(x) := x$ when $V(x) = 0$. Let $\overline{\sigma}$ denote the completion of a $\sigma$-field.

**Theorem 2.1.** *Let $\{S_n\}$ be a transient Markov chain starting at a fixed state, $x_0$. Then $\mathcal{E} \subseteq \overline{\sigma}(\{M(x,y), U(x)\,;\, x,y \in X\})$.*

*Proof.* As in Wilson [9], we imagine running the Markov chain by using infinite stacks under each of the states. The stack under a state $x$ consists of possible successors to $x$ and is generated independently of all other stacks by using the transition probabilities from $x$ repeatedly for independent successors. Once the stacks are generated, the chain moves by moving to the state given at the top of the stack under $x_0$ and removing ("popping") the top state under $x_0$. This is repeated from the current state, and so on. The number of states under $x$ that are eventually popped equals $V(x)$ and the last one is $U(x)$. Let $W(x)$ be the ordered list of states under $x$ that are popped, *excluding* the last one. Write $[W(x)]$ for multi-set of states in $W(x)$, i.e., the unordered list of states (with repetition) in $W(x)$. Note that $\sigma(M(x,y), U(x)\,;\, x,y \in X) = \sigma([W(x)], U(x)\,;\, x \in X)$.

We first claim that if $W(x)$ is re-ordered for $x$ in some finite set of states $A$, then the resulting chain $\{S_n'\}$ starting at $x_0$ will have the same counts $M(x,y)$ and same final exits $U(x)$. It suffices to prove this when $A$ is a singleton. Moreover, if $A$ is not $x_0$, then we may simply begin the chain when it first reaches $A$ and pop the states that are used before then, reducing the situation to $A = \{x_0\}$. Thus, let $A = \{x_0\}$. The transitions of the chain $(S_0, S_1, \ldots)$ describe an Eulerian circuit of a directed multi-graph, $G$. That is, $G$ consists of directed edges $(S_k, S_{k+1})$ connecting vertices $\{S_k\}$ and each vertex has the same number of edges leading to it as leading away from it, except that $x_0$ has one more edge leading away. When $W(x_0)$ is re-ordered, the sequence $(S_0', S_1', \ldots)$ does not leave $G$ (while using each edge at most once) since the number of possible arrivals to a vertex via an edge of $G$ is at most the number of possible departures. Thus, $(S_0', S_1', \ldots)$ traverses a subgraph $G'$ of $G$. If we re-order again to the original order, then this argument shows that the resulting graph covered, $G$, is a subgraph of $G'$. Thus, $G' = G$. Therefore, the final transition counts are the same, as claimed. In addition, the stacks were popped in the same order at all vertices other than $x_0$, so their final exits are unchanged, as is $U(x_0)$.

We next claim that the distribution of $\{S_n\}$ given $[W(x)]$ and $U(x)$ for all $x \in X$ can be represented as follows: Choose randomly and uniformly an ordering $W(x)$ for each $[W(x)]$, independently for each $x \in X$. Then the resulting walk starting from $x_0$ and determined by these stacks has the same law as the Markov chain. To see this, consider the set $B$ of trajectories that correspond to a given collection of $[W(x)]$ and $U(x)$. Let $\{S_n\} \in B$ be one such trajectory. Since re-ordering any finite set of the corresponding $W(x)$ gives a finite permutation of $\{S_n\}$ with the same counts and final exits, $B$ and the conditional Markov chain measure on $B$ are preserved. Therefore the Markov chain measure is preserved under re-ordering every $W(x)$. The only such invariant measure is the one described, so the claim is proved.

Finally, let $C \in \mathcal{E}$. Let $B$ be the set of trajectories that correspond to a given collection of $[W(x)]$ and $U(x)$. Since both $C$ and $B$ are invariant under re-ordering any finite $W(x)$, so is $C \cap B$. In addition, the orderings $W(x)$ are independent given all $[W(x)]$ (and $U(x)$), so the conditional probability of $C$ given $B$ is 0 or 1 by Kolmogorov's 0-1 law. Let $D_0$ be the union of those $B$ for which the conditional probability of $C$ given $B$ is 0 and $D_1$ be the union of the other $B$. Then $P[C \cap D_0] = 0$, so $P[C \triangle D_1] = 0$. Since $D_1 \in \overline{\sigma}\bigl(M(x,y), U(x) \,;\ x, y \in X\bigr)$, the theorem is proved. □

**Corollary 2.1.** *For a transient nearest-neighbor random walk on a tree (with arbitrary transition probabilities), we have $\overline{\mathcal{E}} = \overline{\sigma}(V)$.*

*Proof.* Since a transient random walk on a tree $T$ must tend to some end of $T$, it follows that the pointers $U(x)$ are determined by the occupation field $V$. In view of the preceding theorem, it suffices to show that the transition numbers $M(x,y)$ are also determined by $V$. Write $L_0 = S_0 = x_0$, and for $j \geq 1$ define $L_j = U(L_{j-1})$. The sequence $L = \{L_j \,;\ j \geq 0\}$ is known as the *loop-erasure* of the trajectory $\{S_k \,;\ k \geq 0\}$. Consider the finite tree $T_F = T_F(L_k)$ that is spanned by $L_k$ and all vertices $x$ with $V(x) > 0$ and that can be reached from $x_0$ without visiting $L_k$. The proof will now follow from the following **claim:** *Given a finite walk from $x_0$ to $y$ on a finite tree $T_F$, the edge transition numbers $M_F$ of the walk are determined by the occupation numbers $V_F$ of all vertices except $y$.* The claim is proved by induction on the number $N$ of vertices in $T_F$. The base case $N \leq 2$ is clear. For $N > 2$, the tree $T_F$ has some leaf $z$ that is different from $y$. Let $z_*$ denote the neighbor of $z$. Clearly $M(z, z_*) = V(z)$ and $M(z_*, z) = V(z) - \mathbf{1}_{z = x_0}$. Removing $z$ from the tree and subtracting $V_F(z)$ from $V_F(z_*)$ reduces the problem to a tree with $N-1$ vertices and completes the induction step. To apply the claim to our situation, take $y = L_k$ and observe that for all vertices $w \in T_F(L_k)$ except possibly $L_k$ itself, the occupation number $V(w)$ determined by the infinite random walk path coincides with $V_F(w)$, the occupation number determined by the portion of that path in $T_F(L_k)$. (It is certainly possible that $V(L_k) > V_F(L_k)$, due to excursions of the random walk from $L_k$ to the complement of $T_F$.) □

## 3. A transient birth-and-death chain with finitely many cutpoints

We shall exhibit a birth-and-death chain, i.e., a nearest-neighbor random walk on $\mathbf{N}$, which is transient but has only finitely many cutpoints a.s. We shall use the following basic fact about random walks and electrical networks. Let $r_k > 0$ be given for $k \geq 1$. (Interpret $r_k$ as the resistance of the edge between $k$ and $k+1$.) Consider the birth-and-death chain on $\{1, 2, \ldots, n\}$ where the transition probability from 1 to 2 is 1, and for $k > 1$, the transition probability from $k$ to $k+1$ is $r_{k-1}/(r_{k-1}+r_k)$ and the transition probability from $k$ to $k-1$ is $r_k/(r_{k-1}+r_k)$. Then the probability that the chain reaches $n$ before 1 when starting from $k$ equals $\sum_{j=1}^{k-1} r_j / \sum_{j=1}^{n-1} r_j$. See [2], §§II.1 and IX.2. Of course, this can also be phrased as a standard gambler's ruin calculation. In particular, taking a limit as $n \to \infty$ shows that transience is equivalent to $\sum_{j=1}^{\infty} r_j < \infty$.

**Theorem 3.1.** *Fix $\beta > 1$. Let $r_k > 0$ have the property that $r_k \asymp k^{-1}(\log k)^{-\beta}$ for all $k \geq 2$, where the symbol $\asymp$ means that the ratio of the two sides is bounded above and below by positive constants that do not depend on $k$. Consider the birth-and-death chain on $\mathbf{N} = \{1, 2, \ldots\}$ with transition probability $r_{k-1}/(r_{k-1}+r_k)$ from $k$ to $k+1$ and transition probability $r_k/(r_{k-1}+r_k)$ from $k$ to $k-1$ for all $k \geq 2$. (The transition probability from 1 to 2 is 1.) Then this chain is transient and has only finitely many cutpoints a.s.*

*Proof.* We may assume the chain starts at 1. Since $\sum_k r_k < \infty$, the walk is transient.

Denote $t_k := \sum_{j \geq k} r_j$. The usual gambler's ruin calculation shows that the probability that the walk will have $k$ as a cutpoint is $p_k = r_k/t_k$.

Let $j < k$. Given that $k$ is a cutpoint, let $Q_k(j)$ be the conditional probability that $j$ is a cutpoint. Then $Q_k(j)$ is the probability that a walk starting at $j+1$ visits $k+1$ before visiting $j$, i.e.,

$$(4) \quad Q_k(j) = \frac{r_j}{(t_j - t_{k+1})}.$$

This is also the conditional probability

$$\mathbf{P}[j \text{ is a cutpoint} \mid k \text{ is a cutpoint}, F_{k+1}],$$

where $F_{k+1}$ is any event determined by the future of the walk after it reaches $k+1$ for the first time.

Let $C_{j,k}$ be the set of cutpoints in $(2^j, 2^k]$ and $A_{j,k} := |C_{j,k}|$. Write $a_m := P[A_{m,m+1} > 0]$ and

$$b_m := \min\Big\{\sum_{i=1}^{2^{m-1}} Q_k(k-i) \,;\, k \in (2^m, 2^{m+1}]\Big\}.$$

On the event that $A_{m,m+1} > 0$, let $\ell_m$ be the largest cutpoint in $C_{m,m+1}$. Bound below the expected number of cutpoints in $(2^{m-1}, 2^{m+1}]$ by conditioning on the last cutpoint in $(2^m, 2^{m+1}]$, if there is one:

$$(5) \quad \sum_{j=2^{m-1}+1}^{2^{m+1}} p_j = \mathbf{E}[A_{m-1,m+1}]$$
$$\geq a_m \mathbf{E}[A_{m-1,m+1} \mid A_{m,m+1} > 0]$$
$$= a_m \mathbf{E}\big[\mathbf{E}[A_{m-1,m+1} \mid A_{m,m+1} > 0, \ell_m]\big]$$
$$\geq a_m b_m.$$

Now $t_j \asymp (\log j)^{-\beta+1}$, whence $p_j \asymp (j \log j)^{-1}$ for $j \geq 2$. Furthermore, we have $t_{k-i} - t_{k+1} \asymp i r_k \asymp i r_{k-i}$ for $1 \leq i \leq 2^{m-1}$ and $2^m < k \leq 2^{m+1}$. By (4), this means that $Q_k(k-i) \geq c/i$ for some constant $c > 0$ and $i, k$ in those ranges, which gives in turn that $b_m \geq c'm$ for some constant $c' > 0$. On the other hand, the left-hand side of (5) is at most $c''(\log\log 2^{m+1} - \log\log 2^m) \leq c'''/m$ for some $c'', c''' < \infty$. It follows that $a_m = O(1/m^2)$ is summable, so that there are a.s. only finitely many cutpoints by the Borel-Cantelli lemma. It also follows that with positive probability, there are no cutpoints at all. □

## 4. Concluding remarks

Given a transient Markov chain $\{S_j\}$ with a fixed starting state, it is easy to see that for any $n$, the event $A_n$ that $S_0, S_1 \ldots, S_n$ are all cutpoints has positive probability. Indeed, starting from a trajectory $S_0, S_1, S_2 \ldots$, consider the corresponding loop-erased path $\{L_j\}$ obtained by erasing cycles in the path as they are created. More precisely, $L_0 = x_0$ and $L_j = U(L_{j-1})$ for $j > 0$, where $U(\cdot)$ is the ultimate successor function defined in Section 2. Fix a sequence of vertices $(x_1, \ldots x_n)$ such that the event $B_n = \{(L_0, \ldots, L_n) = (x_0, \ldots, x_n)\}$ has $P(B_n) > 0$. If $B_n$ holds for the trajectory $\{S_j^*\}$, then $x_j = L_j = S_{k_j}$ for some random sequence $\{k_j\}$, and we define a new trajectory $\{S_j^*\}$ with $S_j^* = L_j$ for $j = 0, \ldots, n$ and $S_{n+i}^* = S_{k_n+i}$ for $i > 0$. For this new trajectory $x_0, \ldots, x_n$ are all cutpoints. We conclude that $P(A_n) \geq P(B_n) \prod_{j=1}^n p(x_{j-1}, x_j) > 0$.

We do not know whether every transient Markov chain has an infinite expected number of cutpoints. For any birth-and-death chain, this does hold since (in the notation of the preceding proof) $\sum_{k \geq m} p_k \geq \sum_{k \geq m} r_k/t_m = 1$ for every $m$, whence the series $\sum_k p_k$ diverges.

Another natural question that we cannot answer is whether a *simple* random walk on any transient graph of bounded degree must have infinitely many cutpoints a.s.